\documentclass[prb,twocolumn,showpacs,amsmath,amssymb,superscriptaddress]{revtex4-2}
\usepackage{amssymb}
\usepackage{amsmath}
\usepackage{colordvi}
\usepackage[colorlinks]{hyperref}
\usepackage{amsthm}
\usepackage{subeqnarray}
\usepackage{cases}
\usepackage{enumerate}
\usepackage{graphicx}
\usepackage{epsfig}
\usepackage{epstopdf}
\usepackage{subfigure}
\usepackage{color}
\usepackage{CJKutf8}
\usepackage{tikz}
\usetikzlibrary {arrows.meta,positioning,hobby}

\def\avg(#1){\langle#1\rangle}

\def\be{\begin{equation}}
\def\ee{\end{equation}}
\def\bea{\begin{eqnarray}}
\def\eea{\end{eqnarray}}

\DeclareMathOperator{\coker}{coker}
\DeclareMathOperator{\im}{im}
\DeclareMathOperator{\coim}{coim}

\begin{document}
\title{Generalization to the Natural Gradient Descent}

\author{Shaojun Dong}
\affiliation{Institute of Artificial Intelligence, Hefei Comprehensive National Science Center}

\author{Fengyu Le}
\affiliation{AHU-IAI AI Joint Laboratory, Anhui University}
\affiliation{Institute of Artificial Intelligence, Hefei Comprehensive National Science Center}

\author{Meng Zhang}
\affiliation{CAS Key Laboratory of Quantum Information, University of Science and Technology of China, Hefei 230026, People's Republic of China}
\affiliation{Synergetic Innovation Center of Quantum Information and Quantum Physics, University of Science and Technology of China, Hefei 230026, China}

\author{Si-Jing Tao}
\affiliation{CAS Key Laboratory of Quantum Information, University of Science and Technology of China, Hefei 230026, People's Republic of China}
\affiliation{Synergetic Innovation Center of Quantum Information and Quantum Physics, University of Science and Technology of China, Hefei 230026, China}

\author{Chao Wang}
\email{wang1329@ustc.edu.cn}
\affiliation{Institute of Artificial Intelligence, Hefei Comprehensive National Science Center}

\author{Yong-Jian Han}
\email{smhan@ustc.edu.cn}
\affiliation{CAS Key Laboratory of Quantum Information, University of Science and Technology of China, Hefei 230026, People's Republic of China}
\affiliation{Institute of Artificial Intelligence, Hefei Comprehensive National Science Center}
\affiliation{Synergetic Innovation Center of Quantum Information and Quantum Physics, University of Science and Technology of China, Hefei 230026, China}

\author{Guo-Ping Guo}
\affiliation{CAS Key Laboratory of Quantum Information, University of Science and Technology of China, Hefei 230026, People's Republic of China}
\affiliation{Institute of Artificial Intelligence, Hefei Comprehensive National Science Center}
\affiliation{Synergetic Innovation Center of Quantum Information and Quantum Physics, University of Science and Technology of China, Hefei 230026, China}

\begin{abstract}
Optimization problem, which is aimed at finding the global minimal value of a given cost function, is one of the central problem in science and engineering. Various numerical methods have been proposed to solve this problem, among which the Gradient Descent (GD) method is the most popular one due to its simplicity and efficiency. However, the GD method suffers from two main issues: the local minima and the slow convergence especially near the minima point. The Natural Gradient Descent(NGD), which has been proved as one of the most powerful method for various optimization problems in machine learning, tensor network, variational quantum algorithms and so on, supplies an efficient way to accelerate the convergence. Here, we give a unified method to extend the NGD method to a more general situation which keeps the fast convergence by looking for a more suitable metric through introducing a 'proper' reference Riemannian manifold. Our method generalizes the NDG, and may give more insight of the optimization methods.
\end{abstract}
\maketitle

\section{introduction}

Optimization problems play crucial role in many fields of science, engineering and industry.  Generally, a task can be evaluated by a real-valued function $L({x})$ (named cost function), where ${x}$ are the parameters used to identify the solution. Therefore, the key to solve the problem optimally is searching the parameters $x$ minimizing the function $L({x})$. It is a big challenge to find the global minima of the cost function $L({x})$ when it is complicated.

Various numerical methods have been introduced to solve the optimization problems. Among which, the gradient descent (GD) method(also known as steepest descent method), is a widely used technique for its simplicity and high efficiency. In the GD method, the parameters are updated according to the gradient, reads
\begin{equation}\label{equ::GM}
x_{t+1}=x_{t} - \eta \partial_x L
\end{equation}
where $\eta$ is the searching step at iteration $t$. Generally, numerical methods, which reach global minima by iterative updates of parameters according to local information, suffer two main problems: the local minima and the slow convergence. The stochastic methods are introduced to help the method to jump out the local minima\cite{Adam2014}. On the other hand, the convergence has been improved by more sophisticated methods like Conjugate Gradient(CG) method\cite{CG0,CG1,CG}, Adam method\cite{Adam2014}, Natural Gradient Descent(NGD) method\cite{Amari1998NGD,Martens2020GD}. 

Let's take another perspective to view the gradient descent:  Let $X$ be the parameter space which form a manifold with a flat metric $G_{ij}(x)=\delta_{ij}$. $L(x)$ is viewed as a scalar field defined on $X$. It's easy to see that $-\partial_x L$ is the direction in which we can obtain largest decrease of cost function for fixed small step length, where the distance is defined by the flat metric. Following this idea, the NGD replace the flat metric  on the parameter manifold $X$ by a generally non-flat metric $G$, in order to improve the convergence. The optimization variables update with the iteration
\begin{equation}\label{equ::NG}
{ x}_{t+1}={ x}_{t} - \eta G^{-1}\frac{\partial L({ x})}{\partial { x}}.
\end{equation}
The most common application of NGD is the optimizations  of statistical models, where the Fisher information(FI) matrix \cite{Amari1998NGD,Martens2020GD}is used as the metric. Fisher information(FI) matrix for parameterized probability distribution $p(x,s)$ is defined as
\begin{equation}\label{metric::FI}
G^{FI}_{i,j}(x)=\sum_sp(x,s)\frac{\partial \log p(x,s)}{\partial x^i}\frac{\partial \log p(x,s)}{\partial x^j}
\end{equation}
where $x$ is the parameter and $s$ is the random variable.

NGD has gained more and more attentions in recent years by showing its outstanding performance in convergence on a variety areas of researches such as Variational Quantum Algorithms(VQA)\cite{Wecker2015,Broers2021,Haug2021,Yao2022},Variational Quantum Eigensolver(VQE)\cite{Wierichs2020,Li2017,Gidi2022,Koczor2019}.  The NGD has also been used in solving quantum many-body problems\cite{Liang2021} where the model can be transformed to a statistical one, and the metric is called Fubini-Study metric tensor\cite{Stokes2020,Koczor2019}. In the field of variational Monte Carlo method in neural network\cite{Shi2019,Nagy2019}, NGD (also called Stochastic Reconfiguration\cite{Sorella1998,Park2020,Park2000}) method is also widely used, and the FI matrix is called S matrix. Although the NGD has shown its efficiency in various realm, there comes a question whether the FI matrix is the only choice for the metric in the NGD? Could we find a better one? Moreover, there are many optimization problems which can not be transformed into statistical models naturally. How can we find the proper metric used in NGD in such problems?

In this paper, we generalize natural gradient descent method by proposing a systematic method to find out a class of well-performed metrics on a broader class of cost functions. Firstly, we introduce a reference space that is highly relevant to the cost function, such that a good metric is easy to choose on the reference space. Then the metric on the parameter space used in NGD can be chosen as the induced metric from the good metric on the reference space. Our method is benchmarked on several examples and is proved to converge faster than GD and CG method. The results also showed that out metrics have better performance than the Fisher information matrix even on some statistical models.

The rest of the paper is organized as follows. We first discuss the motivation and effectiveness of our algorithm by capturing the geometry of the optimization landscape in sec.(\ref{section::theory}). Our algorithm is explained in detail in sec.(\ref{section::algorithm}). Then we show some examples in the sec.(\ref{section::Numerical}). In the first example in subsection.(\ref{subsection::LSM}), our algorithm is applied to the least square problem, which is common in the realm of the deep learning and the statistics. In this example, we found at least 4 distinctive metrics for the NGD and all of them work with high efficiency. In the second example, our algorithm is applied to a classical spin model of Heisenberg model. The FI matrix fails in the third example in subsection.(\ref{subsection::eig}), while our method can provide a workable metric of high quality. We will give a summary in sec.(\ref{section::conclusion}).

\section{theory}\label{section::theory}

In this section, we explain our motivation in determining the metric in manifold $X$ by introducing a reference manifold $Y$. We also try to analysis the effectiveness of our method.

\begin{figure}
	\centering
	\begin{tikzpicture}
		\draw (0,0) ellipse [x radius=2, y radius=1];
		\draw[->] (0,0) -- (1.95,0);
		\filldraw (0,0) circle [radius=0.3mm];
		\node[below] at (0,0) {$x$};
		\node[above] at (1,0) {$dx$};
		\node at (2,1.1) {$|dx|=\epsilon$};
		\node at (0,-1.5) {$X$};\
	\end{tikzpicture}
	\caption{The constraint in determining the gradient direction at $x\in X$.}\label{fig:constraint_x}
\end{figure}
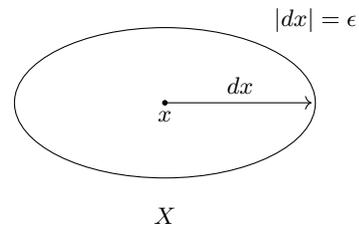

Firstly we follow the line in Ref.\onlinecite{Amari1998NGD} to give a brief introduction to the natural gradient descent method. Consider a cost function $L:X\rightarrow \mathbb R$, where $X$ is the parameters space. For every point $x\in X$, we define a metric $G_X(x)$ at $x$. Starting from some initial parameters $x_{0}\in X$, we search for $x_{min}\in X$ to make $L(x)$ minimal. At parameter $x_t\in X$, the strategy to descent $L(x)$ is to perform a line search in the direction $dx_t$ (which is dependent of the parameter $x_{t}$)
\begin{equation}
	x_{t+1}=x_t+\eta dx_t \label{eqn:evolve_x},
\end{equation}
where $dx_t$ is the direction in which we obtain the maximal descent of $L(x)$ by performing an update with fixed step size $|dx|=(G_{X,ij}dx^idx^j)^{1/2}$. $dx_t$ can be determined by the following fomula
\begin{align}
	dx_t&=-\frac 1\epsilon\underset{dx}\max [L(x_t+dx)] \quad dx\in TX(x_t),|dx|=\epsilon\\
  &=-\frac 1\epsilon\underset{dx}\max [dx^i\frac{\partial L}{\partial x^i}]\quad dx\in TX(x_t),|dx|=\epsilon
\end{align}
where $TX(x_t)$ is the tangent space at $x_t$. The constraint is shown in Fig.\ref{fig:constraint_x}.

 The minimization can be done using Laplacian multiplier method, and the solution $dx_t$ is given by\cite{Amari1998NGD}
\begin{equation}
	 dx_t=-G_{X}^{-1}(x_t)\partial_x L \label{eqn:dir_x}
\end{equation}
up to a positive normalizing factor.

However, $dx_t$ only gives a direction that $L(x)$ decrease fast locally. There's no guarantee that the direction $dx_t$ obtained this way is the optimal at global scale. It's clear to see that at each point $x\in X$, the globally best evolution direction is $\overrightarrow{x,x_{min}}$, where $x_{min}$ is the global minimum. This can be used to judge the global effectiveness of a direction (and the metric). We define a good metric $G_X$ for the cost function $L(x)$ as one that satisfies
\begin{equation}
	dx=-G_X^{-1}\partial_x L(x)\sim \overrightarrow{x,x_{min}}\label{eqn:criterion_y}
\end{equation}

However natural gradient descent method bears an issue that a good metric $G_X$ is often hard to find due to the complexity of $L(x)$. At present, we can only write metric for some special types of problem such as statistical model optimization and quantum eigensolver, based on case-specific study. 

In this work, we propose a systematic method to provide metric of good quality for a class of cost functions: one that can be re-written as $L(x)=\bar L(f(x))$, where $f:X\rightarrow Y$ is a map from parameter manifold $X$ to a reference manifold $Y$ (which is also a Riemannian manifold), such that the following conditions are satisfied:
\begin{enumerate}
	\item $\bar L(y)$ is a relatively simple function such as a polynomial function
 or a function whose well-performed metric is known. Therefore we can find a  good metric $G_Y$ in $Y$ manifold for the cost function $\bar L(y)$.
	\item $f$ is a relatively surjective map especially when the cost function $L(x)$ is small.
	\item Let $x_{min}\in X$ be the minimum of $L$, and $y_{min}\in X$ be the minimum of $\bar L$. Then $f(x_{min})$ is close to $y_{min}$. 
\end{enumerate}

Such class of cost functions include the ansatz-type problems, where we often need to optimize a simple cost function $\bar L(p)$ over a space $P$ of physical meaningful quantities. However, due to high-dimensionality of $p\in P$, the elements in $P$ is simplified by some ansatz $f:A\rightarrow P$, which has strong expressibility when $\bar L(p)$ is small. Due to high complexity of $f$, a direct choice of metric for function $\bar L(f(a))$ is often of poor quality. For ansatz-type problems, we may use $P$ as reference space and $A$ as parameter space. The requirements above are commonly satisfied.

Given a metric $G_Y$ in $Y$ manifold, we'll show that define $G^f_X$ by the metric induced by $f$
\begin{equation}
	G^f_{X,ij}(x)=\frac{\partial y^\alpha}{\partial x^i}G_{Y,\alpha\beta}(f(x))\frac{\partial y^\beta}{\partial x^j} \label{eqn:GX}
\end{equation}
is a good metric for the cost function L(x).

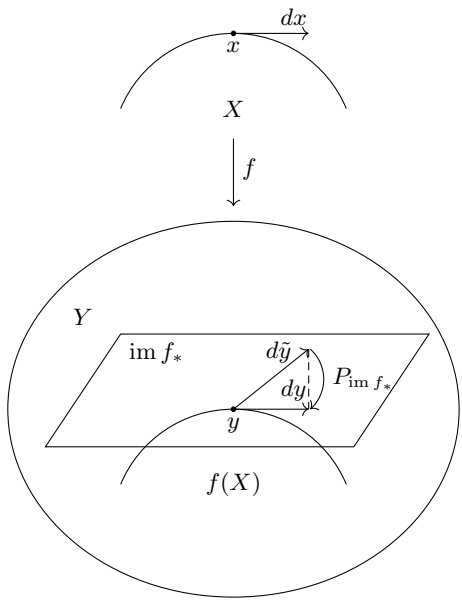
\begin{figure}
	\centering
	\begin{tikzpicture}
		\draw [use Hobby shortcut] (-1.5,-1) .. (0,0) .. (1.5,-1);
		\draw[->] (0,0) -- (1,0);
		\filldraw (0,0) circle [radius=0.3mm];
		\node[below] at (0,0) {$x$};
		\node[above] at (0.8,0) {$dx$};
		\node at (0,-1) {$X$};
		\draw[->] (0,-1.4) -- (0,-2.3);
		\node[right] at (0,-1.8) {$f$};
		\begin{scope}[yshift=-5cm]
		\draw (0,0) ellipse [x radius=3, y radius=2.5];
		\draw [use Hobby shortcut] (-1.5,-1) .. (0,0) .. (1.5,-1);
		\draw[->] (0,0) -- (1,0);
		\draw[->] (0,0) -- (1,0.8);
		\draw[densely dashed] (1,0.8) -- (1,0);
		\draw (-2.5,-0.5) -- (1.6,-0.5) -- (2.6,1) -- (-1.5,1) -- cycle;
		\draw[->,use Hobby shortcut] (1.02,0.8) .. (1.2,0.4) .. (1.02,0);
		\filldraw (0,0) circle [radius=0.3mm];
		\node[below] at (0,0) {$y$};
		\node[above] at (0.8,0) {$dy$};
		\node[right] at (1.2,0.4) {$P_{\im f_*}$};
		\node[below right] at (-1.5,1) {$\im f_*$};
		\node at (0,-1) {$f(X)$};
		\node[above] at (-2,1) {$Y$};
		\node[above] at (0.6,0.5) {$d\tilde y$};
		\end{scope}
	\end{tikzpicture}
	\caption{The relation between $dy=f_*(dx)$ and $d\tilde y =G_Y^{-1}\partial_y\bar L$.}\label{fig:relation_dy_simp}
\end{figure}

At fixed point $x\in X$ and $y=f(x)\in Y$, $f$ induce a linear map $f_* :dx\mapsto dy$ defined by
\begin{equation}
	dy^\alpha=\frac{\partial y^\alpha}{\partial x^i}dx^i ,
\end{equation}
 which maps a tangent vector $dx$ at $x$ to a tangent vector $dy$ at $y$. Given a line search direction $dx$ determined by Eq.\ref{eqn:dir_x}, as proved in Append.\ref{section::app1}, it is mapped to an update in reference manifold in the direction $dy=f_* (dx)\propto P_{\im f_*}G_Y^{-1}\partial_y\bar L$, where $P_{\im f_*}$ is the projection operator from the tangent space at $y$ to the image of $f_*$ with respect to the metric $G_Y$. The process is shown in Fig.\ref{fig:relation_dy_simp}.

It's easy to see that $G^f_X$ is a good metric in $X$ manifold iff
\begin{equation}
	dy=-P_{\im f_*}G_Y^{-1}\partial_y \bar L(y)\sim  \overrightarrow{y,f(x_{min})}\label{eqn:criterion_y2}
\end{equation}
This is likely true if the three conditions we required are satisfied, because
\begin{enumerate}
	\item If $G_Y$ is a good metric in $Y$ manifold, then $-G_Y^{-1}\partial_y \bar L(y)\sim  \overrightarrow{y,y_{min}}$.
	\item If $f$ is a relatively surjective map especially when the cost function is small, then $-P_{\im f_*}G_Y^{-1}\partial_y \bar L(y)\sim -G_Y^{-1}\partial_y \bar L(y)$.
	\item If $f(x_{min})$ is close to $y_{min}$, then $\overrightarrow{y,f(x_{min})}\sim\overrightarrow{y,y_{min}}$. 
\end{enumerate}

The benefits of this process is that we don't need to find the metric in the parameter manifold directly, whose quality is severely impacted by the complexity of original cost function $L$. In stead, the complexity is included in the transformation $f$,  and has relatively little impact on the quality of the final metric we choose.

\section{Algorithm}\label{section::algorithm}
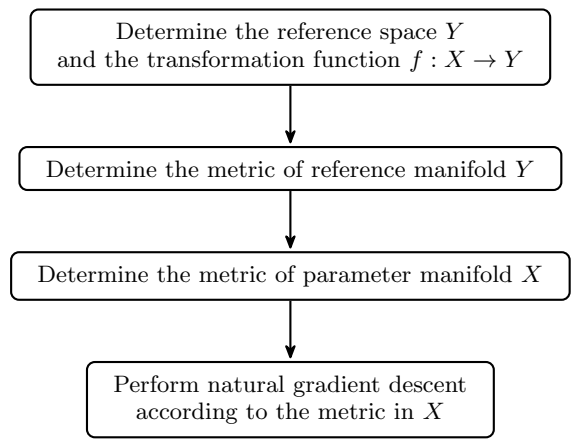
\begin{figure}[htb!]
	\centering
	\begin{tikzpicture}[every node/.style={align=center,rounded corners=3pt,inner xsep=10pt,inner ysep=5pt,draw},->,>={Stealth[round]},shorten >=1pt,thick,node distance=0.8cm]
		\node (1) at (0,0) {Determine the reference space $Y$\\ and the transformation function  $f:X\rightarrow Y$};
		\node (2) [below= of 1] {Determine the metric of reference manifold $Y$};
		\node (3) [below= of 2] {Determine the metric of parameter manifold $X$};
		\node (4) [below= of 3] {Perform natural gradient descent\\ according to the metric in $X$};
		\path (1) edge (2);
		\path (2) edge (3);
		\path (3) edge (4);
	\end{tikzpicture}
	\caption{The scheme of the algorithm.}\label{fig:algorithm}
\end{figure}

In this section, we explain the steps of our algorithm in detail. The algorithm is shown schematically in Fig. \ref{fig:algorithm}.

The first and critical step of our method is to determine the reference manifold $Y$ and the transformation $f:X\rightarrow Y$, such that $L(x)=\bar L(f(x))$, where $\bar L:Y\mapsto \mathbb R$ is a relatively simple function. This calls for experience and understanding of the problem. As we have stated, our method is very suited for the ansatz-type cost function.

Once this is done, we need to determine the metric $G_Y$ in the reference manifold. The most simple choice of $G_Y$ is the Euclidean metric
\begin{equation}
G^{I}_{Y,ij}=\delta_{i,j}
\end{equation}
Another metric taken into consideration in the paper is the Hessian matrix,
\begin{equation}
G^{H}_{Y,ij}=\frac{\partial^2 \bar L({ Y})}{\partial Y_i \partial Y_j}+\epsilon\cdot\delta_{i,j} .
\end{equation}
where the $\epsilon$ is a real number to the ensure the positive-definiteness of the metric $G^{H}_Y$. In our numerical test, $\epsilon=|\epsilon_{H}|+0.1$, where $\epsilon_{H}$ is the minimal eigen-value of the Hessian matrix. The roughly value of $\epsilon_{H}$ can be easily obtained by other methods such as gradient descent.

Once we have a metric for the reference manifold, we obtain the metric $G_X$ in parameter manifold through the conversion of coordinates by Eq.(\ref{eqn:GX}). Then $dx=G^{-1}\frac{\partial L({ x})}{\partial { x}}$  is calculated by solving a linear equation
\begin{equation}\label{equ::NGD_direction}
G \cdot dx = \frac{\partial L({ x})}{\partial { x}} .
\end{equation}
It could be solve easily by the CG method, where the computation cost is only the operation of matrix times vector.

Finally we can determine the update of the parameter by line search in this direction using Eq.(\ref{eqn:evolve_x}).

\section{ Numerical Experiments} \label{section::Numerical}

In this section, we will give examples to clarify our method in capturing the geometry of the optimization landscape. The first example we are going to show is about a least square method(LSM). The cost function of LSM is widely used for the trainning of the deeping learning ansatz. In this example we could find at least 4 metrics for NGD from 3 reference manifolds and all of them have out-standing convergence speed comparing to the CG method and the gradient descent. The second example is the classical spin model. And then follows by a model where NG using fisher information matrix has no advantage over gradient descent, but our method give another metric with excellent results.

We have adopted the line search in the optimizations with the searching directions being normalized.

\subsection{ Least square model} \label{subsection::LSM}
The first example we are going to show is a least square model.

%

As an example, we study a one-dimensional quantum anti-ferromagnetic Heisenberg model\cite{Altlandbook}. The Hamiltonian reads
\begin{equation}
	H=\sum_i \sigma_i^x\sigma_i^x + \sigma_i^y\sigma_i^y + \sigma_i^z\sigma_i^z
\end{equation}
where $\sigma_i^\alpha (\alpha=x,y,z)$ is Pauli matrices. Since the Hamiltonian is real, the ground state of this model can be found in real Hilbert space.

 Suppose that the system is in its ground state which is a pure state. We can measure the reduced density matrices of all pair of nearby sites. However, the state is not ideally isolated form the environment, and the state is actually mixed with some noise. Therefore the measured results are actually reduced density matrices corresponding to some mixed states. Now from the measured data, we use the following optimization method to rebuild the noise-free data. This method could be used for the purifications of the noisy quantum circuits.

Suppose the system is noise-free, the pure ground state reads
\begin{equation}\label{equ::MPS}
|\Psi\rangle=\frac{1}{Z}\sum_{\bf s} W_{\bf s}|{\bf s}\rangle ,
\end{equation}
where $|{\bf s}\rangle=|s_1,s_2,\cdots,s_L\rangle$ are the physical bases of the Hilbert space, $W_{\bf s}$ are coefficients the and $Z^2=\sum_{\bf s} |W_{\bf s}|^2$ is the  normalization coefficient. And the reduced density matrix for site $i$ and site $i+1$ is
\begin{equation}
	(D_{i,i+1})_{s_is_{i+1}s'_is'_{i+1}}=\frac{\sum_{s_1,\dots,s_{i-1},s_{i+2},s_L}W_{\bf s}W_{\tilde{\bf s}}}{Z^2}
\end{equation}
where $\tilde{\bf s}=s_1,\dots,s'_i,s'_{i+1},\dots,s_n$.

Let $\widetilde{D}_{i,i+1}$ be the measured reduced density matrices for site $i$ and site $i+1$, we search the  $D_{i,i+1}({ x})$ to  minimize the cost function of
\begin{equation}\label{equ::LSM}
L({ x})=\frac{1}{L} \sum_{i=1}^{L-1} (D_{i,i+1}({ x}) - \widetilde{D}_{i,i+1})^2
\end{equation}

\begin{figure} [htb!]
		\begin{center}
		\includegraphics[width=0.45\textwidth]{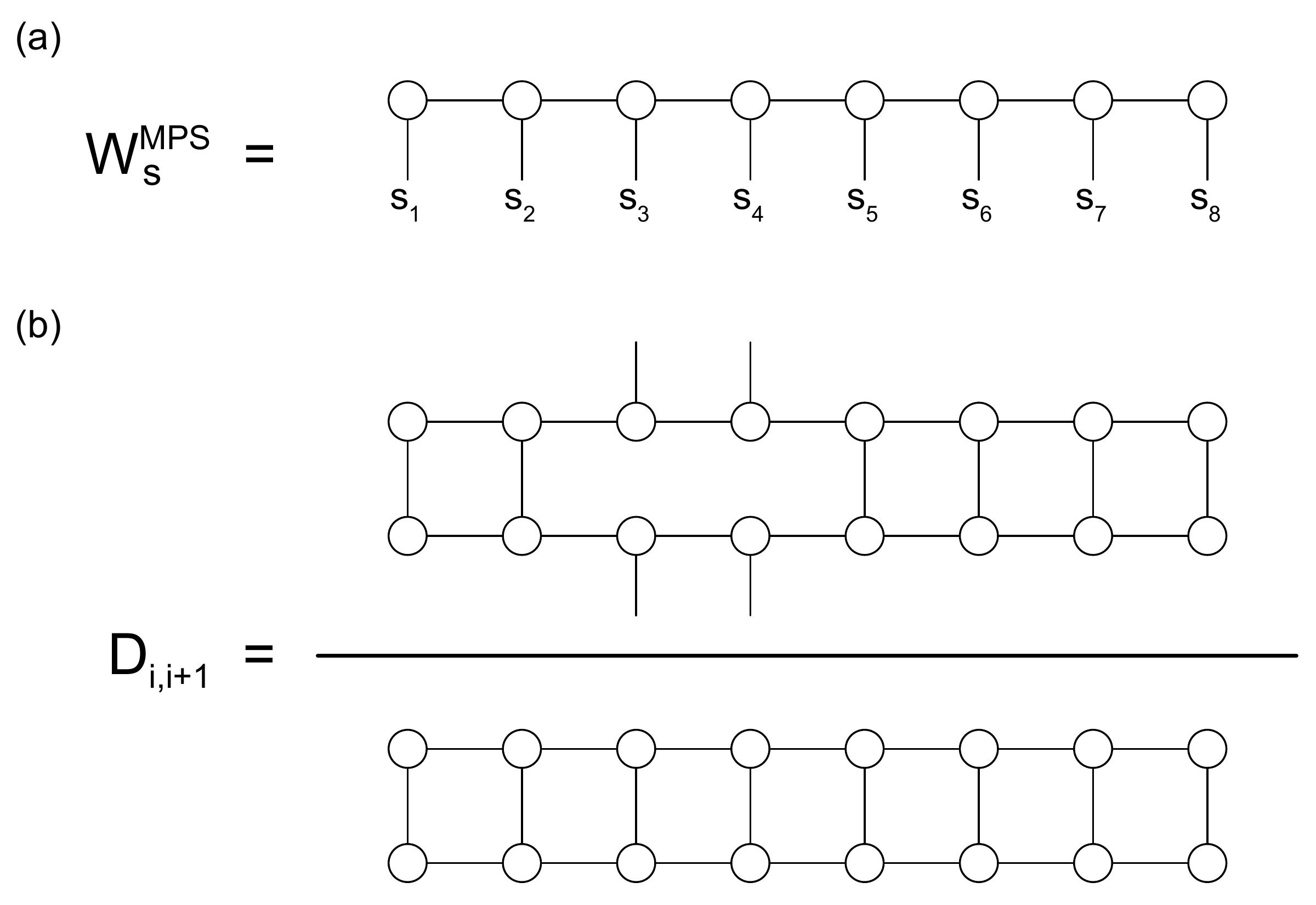}
		\caption{(a) The graphical representations of MPS.  (b) The reduced density matrix $D_{i,i+1}$ for MPS, here $i=3$}\label{pic::MPS}
		\end{center}
\end{figure}

However, the dimension of the Hilbert space is $2^L$. It grows exponentially with the system size. When the system size $L>30$, we can not store and manipulate the exact quantum state in the computer. Here we approximate the quantum state by matrix product state(MPS)\cite{Garcia2006,SCHOLLWOCK201196,Roman2014IntroTNS,Verstraete2008,Orus19}, which is an ansatz to express quantum states with low entanglement. The number of the parameters in the MPS is $\sim 2LD^2$ , which grows polynomially with the system size $L$, where $D$ is the bond dimension of the tensors in the MPS that affection the accuracy of this approximation. In MPS representation, the amplitude $W^{\rm MPS}_{\bf s}$ is of the form
\begin{equation}
	W^{\rm MPS}_{\bf s}=\sum_{a_1,\dots,a_{L-1}}A^{[1]}_{s_1a_1}A^{[2]}_{s_2a_1a_2}\cdots A^{[L-1]}_{s_{L-1}a_{L-2}a_{L-1}}A^{[L]}_{s_La_{L-1}}
\end{equation}
where $A$s are tensors whose elements are free variables(denoted by $ x$), and $a_i$ is summed form 1 to $D$. This can be expressed by diagrams \cite{Orus19} in Fig.\ref{pic::MPS}(a). And the graph representations of the density matrix $D_{i,i+1}({ x})$ is shown in Fig.\ref{pic::MPS}(b).

We apply our optimization method in this problem. There're many available choices of reference space in this problem. As the first choice, we set $Y_i=D_{i,i+1}({ x})$ as the reference space. The cost function on reference space is $L({ Y})=\frac{1}{L} \sum_{i=1}^{L-1} (Y_i - \widetilde{D}_{i,i+1})^2$,
which has a quadratic form. Note that here $Y_i$ are the independent variables in the reference manifold. The identity matrix could be a good choice for the metric in the reference manifold. Then from the Eq.(\ref{eqn:GX}), the metric in the cost function is
\begin{equation}\label{metric::LSM_D}
G^D_{i,j}=\sum_{k} \frac{\partial D_{k,k+1}}{\partial x_i}\frac{\partial D_{k,k+1}}{\partial x_j},
\end{equation}

The second reference space we selected is the Hilbert space, where the independent variables  are $Y_{\bf s}=\frac{W^{\rm MPS}(S)}{Z}$. The cost function on reference space becomes a quartic form: $L({ Y})=\frac{1}{L} \sum_{i}^{L-1} (\sum_{s_1,s_2,\cdots,s_{i-1},s_{i+2},\cdots}Y_{\bf s}Y_{\tilde{\bf s}} - \widetilde{D}_{i,i+1})^2$. For comparison, we test two metrics in this reference manifold, which are the identity matrix and the Hessian matrix. The metrics in the parameter manifold are
\begin{eqnarray}
G^I_{i,j}({ X})&=&\sum_{s} \frac{\partial Y_s}{\partial x_i}\frac{\partial Y_s}{\partial x_j} \label{metric::LSM_IH1} \\
G^H_{i,j}({ X})&=&\sum_{s,s'} \frac{\partial Y_s}{\partial x_i}(\frac{\partial^2 L({ Y})}{\partial Y_s \partial Y_{s'}}+\epsilon\delta_{s,s'})\frac{\partial Y_{s'}}{\partial x_j} \label{metric::LSM_IH2}
\end{eqnarray}

The third reference manifold is the  space of the tensor network state, that is $Y_{\bf s}= W^{\rm MPS}({\bf s})$. As the discussion before, we select the identity matrix as the metric in the reference manifold. The metric in original cost function is
\begin{equation}\label{metric::LSM_S}
G^W_{i,j}=\sum_{s} \frac{\partial W^{\rm MPS}(s)}{\partial x_i}\frac{\partial W^{\rm MPS}(s)}{\partial x_j},
\end{equation}

On a system of size $L=40$, we mimic noise-impacted reduced density matrix (the input data) by the reduced density matrices from a pure ground state of the Heisenberg model, with random numbers in $ (-0.1,0.1)$ added to every elements.  Then we rebuild the pure-state reduced density matrices by a MPS of $D=5$. We optimize the cost function Eq.(\ref{equ::LSM}) with Gradient descent(Eq.(\ref{equ::GM})), CG method and the Natural gradient descent method(Eq.(\ref{equ::NG})) with the metrics define from Eq.(\ref{metric::LSM_D}), Eq.(\ref{metric::LSM_IH1}), Eq.(\ref{metric::LSM_IH2}) and Eq.(\ref{metric::LSM_S}) respectively.

\begin{figure} [tbp]
		\begin{center}
		\includegraphics[width=0.4\textwidth]{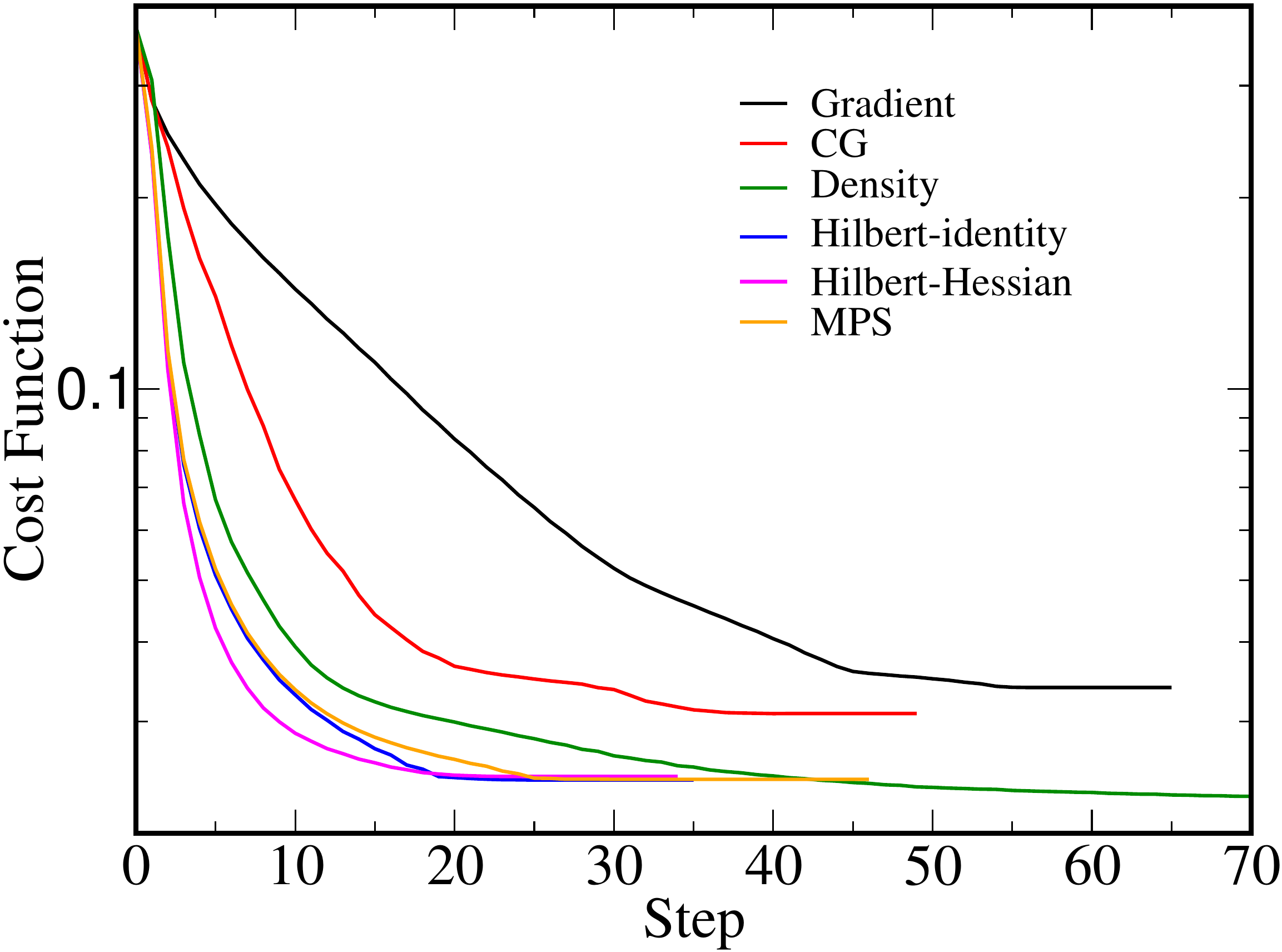}
		\caption{The cost function in Eq.(\ref{equ::LSM}) vs optimization step. We have plot the process of the gradient descent(black line), the CG method(red line) and the NGD with the metrics we find. The label ``Density" with green color is the metric Eq.(\ref{metric::LSM_D}) from the reference manifold of the density matrix. The one ``Hilbert-identity" and the ``Hilbert-Hessian" are the metrics of Eq.(\ref{metric::LSM_IH1},\ref{metric::LSM_IH2}) which are extracted form the Hilbert space of the quantum system. The label ``MPS" is the metric of Eq.(\ref{metric::LSM_S}), which is from the tensor network state. All of our methods are better then the Gradient descent and the CG method both in the convergence speed and the final value.}\label{pic::LSM}
		\end{center}
\end{figure}

The numerical results are shown in Fig.\ref{pic::LSM}, where all the optimizations are terminated when trapped into the local minimums. As shown in Fig.\ref{pic::LSM}, the CG method has better performance than the gradient descent. While all the NGD with the metrics we found are much better than the CG method and the gradient descent both in convergence speed and the final values of the cost function. It seems that the NGD methods can more easily  escape from the local minimum than others. We list all the reference manifolds we use in Tab.\ref{tab:LSM}.

\begin{table*} [t]
\caption{ We list the reference manifolds we used for the cost function Eq.(\ref{equ::LSM}).}

\begin{tabular}{  c c c }
\hline
reference manifold & transformation & cost function \\
\hline
\hline
reduced density matrices   & $Y_i=D_{i,i+1}(x)$ &  quadratic form \\
Hilbert space   & $Y_s=\frac{W(S)}{Z}$ &  quartic form \\
tensor network space   & $Y_s=W(S)$ &  elementary function  \\
\hline
\hline
\end{tabular}
\label{tab:LSM}
 \end{table*}

\subsection{Classical spin model}\label{subsection::spin}
The next numerical experiment is the classical spin model\cite{Onsager1944,Wu1982,Joyce1967}. The classical spin model is important for two reasons. On one hand the classical spin model can be used to describe (thermal) phase transitions at high temperature, where the quantum fluctuations are not important. On the other hand, the classical spin model is the quantum spin model in large S limit\cite{Yang1952,Wu1982,Kosterlitz_1973}. It may give us information of quantum results in some senses.  In this example, we only focus on the comparisons of the performance of  the optimization algorithms.

We test our method on 2D classical Heisenberg model. The model reads,
\begin{equation}\label{equ::spin}
L( {\vec{S}_i})=\frac{1}{N} \sum_{\langle i,j\rangle} \frac{\vec{S}_i}{|S_i|}\cdot\frac{\vec{S}_j}{|S_j|} ,
\end{equation}
where $\vec{S}_i=(S_{i,x},S_{i,y},S_{i,z})$ is a classical spin defined at the $i$-th node, and $|S_i|=\sqrt{S_{i,x}^2 + S_{i,y}^2 + S_{i,z}^2}$ is the length of the spin. The $\sum_{\langle i,j\rangle}$ denotes the summations over all the nearest-neighbor pairs.

To apply our optimization method, we choose $\vec{Y}_i= \frac{\vec{S}_i}{|S_i|}$  as the reference variables and the Eq.(\ref{equ::spin}) become $L( {\vec{Y}_i})=\frac{1}{N} \sum_{\langle i,j\rangle} \vec{Y}_i\cdot\vec{Y}_j $, which is quadratic form. The Hessian matrix would be a good choice as the metric, that is
\begin{eqnarray}\label{matrix::Spin}
G^{spin}_{i,j}({ X})&=&\sum_{s,s'} \frac{\partial Y_s}{\partial x_i}(\frac{\partial^2 L({ Y})}{\partial Y_s \partial Y_{s'}}+\epsilon\delta_{s,s'})\frac{\partial Y_{s'}}{\partial x_j}
\end{eqnarray}

\begin{figure} [tbp]
		\begin{center}
		\includegraphics[width=0.4\textwidth]{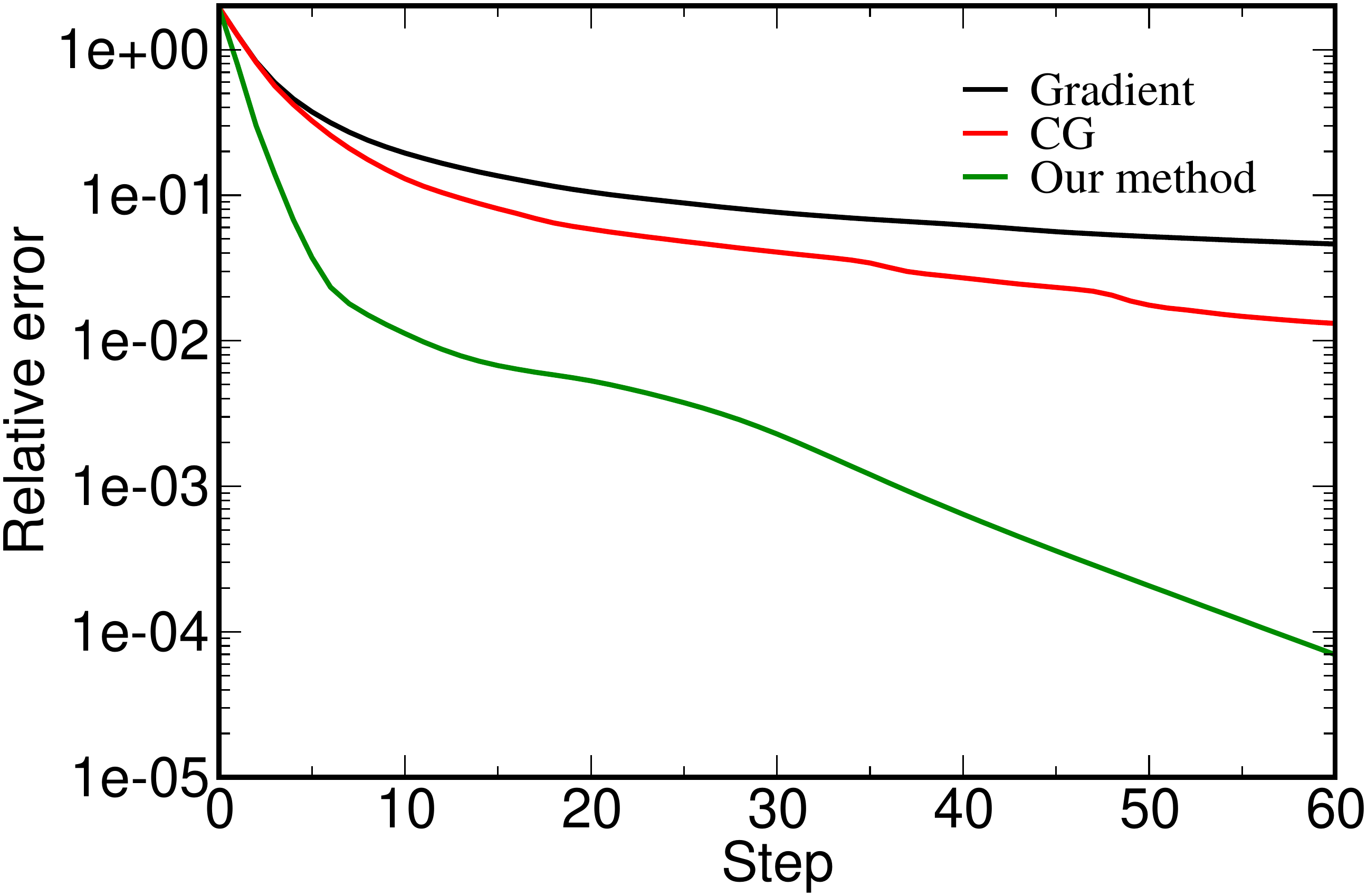}
		\caption{The cost function in Eq.(\ref{equ::spin}) vs optimization step. The relative error is defined as $|(L(\vec{S}_i)-L_{min})/L_{min}|$, where $L_{min}$ is the exact ground state(antiferromagnetic state) energe of Eq.(\ref{equ::Hamiltonian}). Three methods are show for the comparison. They are gradient descent of black line, CG method of red line and the NGD with the matrix we found in the green line. Our method is better then the Gradient descent and the CG method in both the convergence speed and the final  value. }\label{pic::Spin}
		\end{center}
\end{figure}

The numerical results are shown in Fig.\ref{pic::Spin}, where we have tested the Heisenberg model of size $50\times50$ with the Gradient descent, CG method and the NGD with Eq.(\ref{matrix::Spin}) as the metric. Our method has better performance both in convergence speed and final accuracy.

\subsection{Minimal eigen-value of a given matrix}\label{subsection::eig}
Next we are going to show an example in which the well-known FI matrix exists but has no advantage over gradient descent. Fortunately we could find another metric instead through our strategy.

Given a real random matrix(not positive definite) $H$ , we use the optimization methods to find its minimal eigen-value. The cost function is defined as
\begin{equation}\label{equ::Hamiltonian}
L( { x})=\frac{\sum_{i,j} x_iH_{i,j}x_j}{\sum_i x_i^2}
\end{equation}

The FI is calculated through Eq.(\ref{metric::FI}), where the probability in this model is $\rho_i=x_i/Z$ with $Z=\sqrt{\sum_i x_i^2}$. It can be easily proof that the NGD with FI as the metric is identical to the gradient descent in this example.

On the other hand, our strategy can give a metric that work much better that the gradient descent. Here we selected the normalized vector $Y_i=x_i/Z$ as the reference manifold. The cost function on reference space becomes a quadratic function, whose Hessian matrix could be used as the metric. Performing the Eq.(\ref{eqn:GX}),  we have
\begin{equation}\label{matrix::Hamiltonian}
G^{eig}_{i,j}=\sum_{k,l} \frac{\partial Y_k}{\partial x_i}(H_{k,l}+\epsilon\delta_{k,l})\frac{\partial Y_l}{\partial x_j}
\end{equation}

\begin{figure} [tbp]
		\begin{center}
		\includegraphics[width=0.4\textwidth]{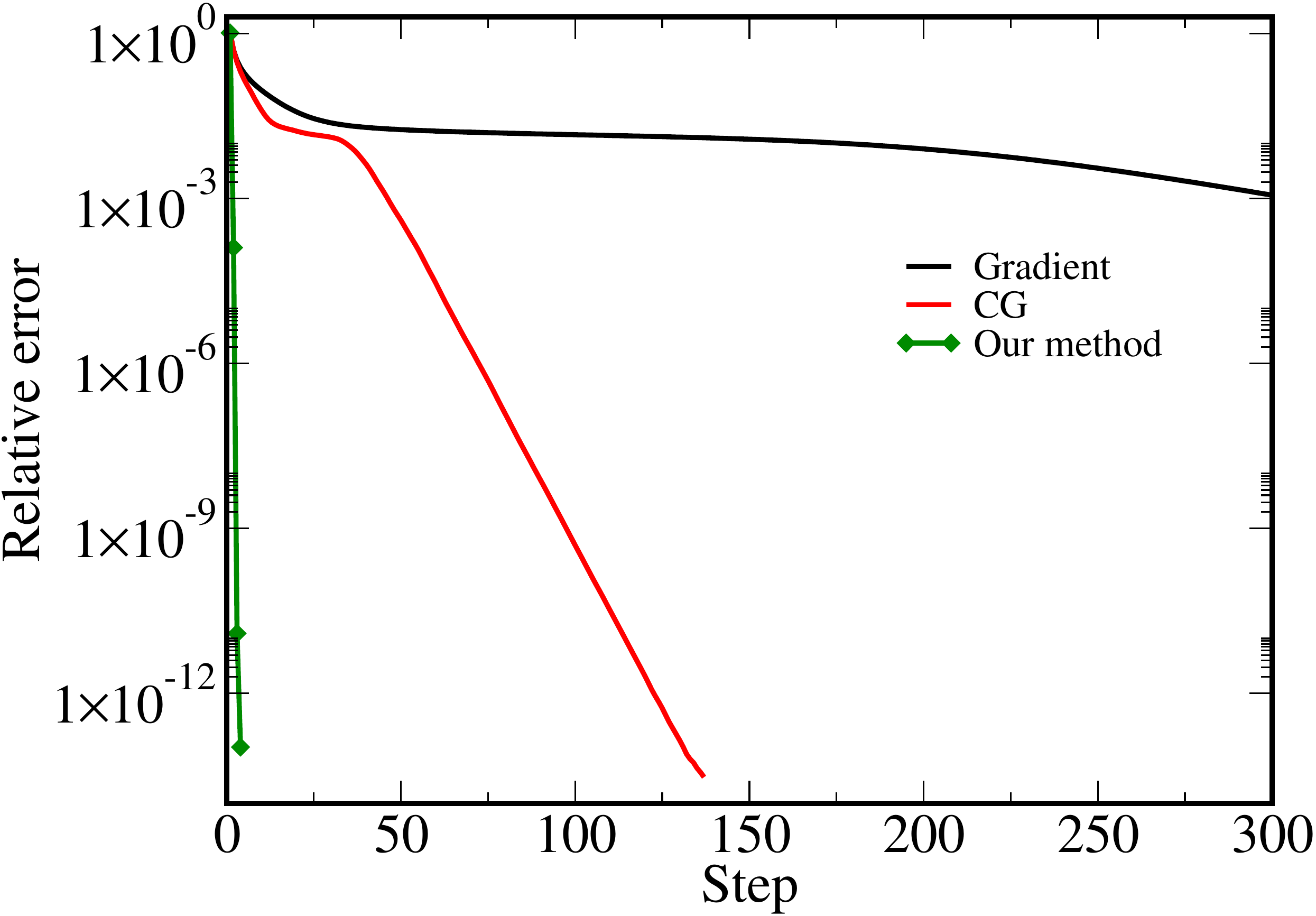}
		\caption{The cost function in Eq.(\ref{equ::Hamiltonian}) vs optimization step. The relative error is defined as $|(L(x)-L_{min})/L_{min}|$, where $L_{min}$ is the exact minimal eigen-value of Eq.(\ref{equ::Hamiltonian}). In the example the search direction from the NGD with FI matrix Eq.(\ref{metric::FI})  are exact the same as that from gradient descent. Three methods are show for the comparison. They are gradient descent in black color, CG method in red line and the NGD with the matrix we found in the green. Here we show that the FI matrix may not be the best choice as the metric for the NGD and our method could find another one to take its place. }\label{pic::eig}
		\end{center}
\end{figure}

The numerical results are shown in Fig.\ref{pic::eig}, where we have tested a $1000\times1000$ random matrix with the Gradient descent(NGD with FI as metric), CG method and the NGD with Eq.(\ref{matrix::Hamiltonian}) as the metric. In this test the CG method is much better than the gradient descent and reached the accuracy of $10^{-12}$ after 125 steps. Our method could reach high accuracy with much fewer steps.

\section{ conclusion} \label{section::conclusion}

We  put forward a  strategy to determine the the geometry of the optimization landscape  used in the natural gradient descent method by looking for a suitable reference manifold to simplify the cost function.

The critical step of our method is to determine the reference manifold $Y$. This can be done by  looking for a  transformation $f:X\rightarrow Y$ as the discussion in sec.(\ref{section::theory}). For the ansatz-type problems, the most simple way to do so is to use the physical space as the reference manifold. We have shown how to do this in several examples. In our examples, more than one metrics  have been found. And all the metrics we found have out-standing performance in the natural gradient descent method even in the case where the Fisher information matrix fails. This proves that our method is more universal than traditional NGD using only Fisher information matrix.

We remark that the purpose of this work is to give more insight into the NGD. We have learnt that the metrics of high quality can be determined by the cost function and  are not unique. This will give us much more freedom  in choosing metrics in the applications of NGD.

\section{Acknowledgement:}

This work was supported by National Natural Science Foundation of China (Nos. 12104433)

\appendix
\section{Update in reference space}\label{section::app1}

In this section, we make connection between $dy=f_*(dx)$ and the gradient-based optimal direction $d\tilde y =G_Y^{-1}\partial_y\bar L$ in reference manifold $Y$.

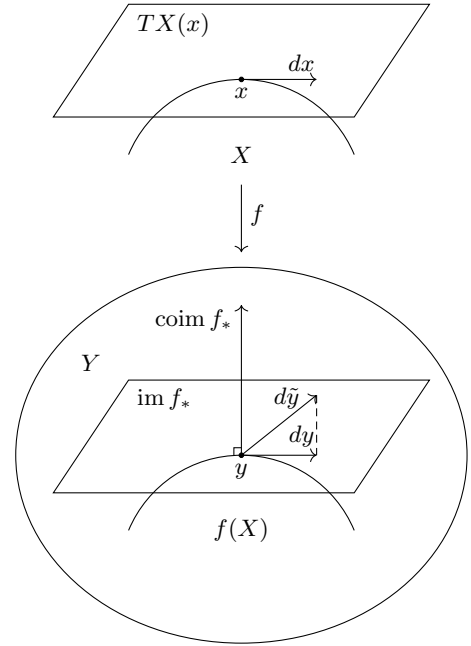
\begin{figure}[htb!]
	\centering
	\begin{tikzpicture}
		\draw [use Hobby shortcut] (-1.5,-1) .. (0,0) .. (1.5,-1);
		\draw[->] (0,0) -- (1,0);
		\draw (-2.5,-0.5) -- (1.5,-0.5) -- (2.5,1) -- (-1.5,1) -- cycle;
		\filldraw (0,0) circle [radius=0.3mm];
		\node[below] at (0,0) {$x$};
		\node[above] at (0.8,0) {$dx$};
		\node[below right] at (-1.5,1) {$TX(x)$};
		\node at (0,-1) {$X$};
		\draw[->] (0,-1.4) -- (0,-2.3);
		\node[right] at (0,-1.8) {$f$};
		\begin{scope}[yshift=-5cm]
		\draw (0,0) ellipse [x radius=3, y radius=2.5];
		\draw [use Hobby shortcut] (-1.5,-1) .. (0,0) .. (1.5,-1);
		\draw[->] (0,0) -- (0,2);
		\draw[->] (0,0) -- (1,0);
		\draw[->] (0,0) -- (1,0.8);
		\draw (0,0.1) -- (-0.1,0.1) -- (-0.1,0);
		\draw[densely dashed] (1,0.8) -- (1,0);
		\draw (-2.5,-0.5) -- (1.5,-0.5) -- (2.5,1) -- (-1.5,1) -- cycle;
		\filldraw (0,0) circle [radius=0.3mm];
		\node[below] at (0,0) {$y$};
		\node[above] at (0.8,0) {$dy$};
		\node[below right] at (-1.5,1) {$\im f_*$};
		\node at (0,-1) {$f(X)$};
		\node[above] at (-2,1) {$Y$};
		\node[left] at (0,1.8) {$\coim f_*$};
		\node[above] at (0.6,0.5) {$d\tilde y$};
		\end{scope}
	\end{tikzpicture}
	\caption{The relation between $dy=f_*(dx)$ and $d\tilde y =G_Y^{-1}\partial_y\bar L$, in detail}\label{fig:relation_dy}
\end{figure}

 We can deduce that
\begin{align}
	dy&=f_*(dx)\\
	&\propto -f_*(\underset{dx}\max [L(x+dx)]) \quad dx\in TX(x),|dx|=\epsilon\\
	&\propto-\underset{dy}\max [\bar L(y+dy)] \quad dy\in \im f_*,|dy|=\epsilon\label{eqn:pre_dir_y}
\end{align}
We may decompose the tangent space at $y$ into $TY(y)=\im f_*\oplus \coker f_*$, where $\im f_*$ is the image of $f_*$, and $\coker f_*$ is the orthogonal complement of $\im f_*$ with respect to the inner product $G_Y(y)$. Let $P_{\im f_*}$ and $P_{\coim f_*}$ be the projection operators into these two subspaces. Similar to Eq.\ref{eqn:dir_x}, Eq.\ref{eqn:pre_dir_y} leads to the equation
\begin{equation}
	P_{\im f_*}G_Ydy\propto P_{\im f_*}\partial_y\bar L
\end{equation}
where $dy\in \im f_*$.
Using $P_{\im f_*}G_YP_{\coim f_*}=0$, it's easy to see that
\begin{equation}
	dy\propto-P_{\im f_*}G_Y^{-1}\partial_y\bar L
\end{equation}
is a solution. Assuming that $P_{\im f_*}G_Y(y)P_{\im f_*}$ has full-rank, it's easy to see that the only solution is
\begin{equation}
	dy=-P_{\im f_*}d\tilde y
\end{equation}
up to a positive normalizing factor, where $d\tilde y =G_Y^{-1}\partial_y\bar L$. The whole process is shown in Fig.\ref{fig:relation_dy}.

\bibliography{ NG.bib}

\end{document}